\documentclass{amsart}
\usepackage{microtype}
\usepackage{amssymb}
\usepackage{mathtools}
\usepackage{tikz-cd}
\usepackage{mathbbol} 
\usepackage[shortcuts]{extdash} 

\usepackage[
	backend=biber,
	style=alphabetic,
	backref=true,
	url=false,
	doi=false,
	isbn=false,
	eprint=false]{biblatex}

\renewbibmacro{in:}{}  
\AtEveryBibitem{\clearfield{pages}}  

\DeclareFieldFormat{title}{\myhref{\mkbibemph{#1}}}
\DeclareFieldFormat
[article,inbook,incollection,inproceedings,patent,thesis,unpublished]
{title}{\myhref{\mkbibquote{#1\isdot}}}

\newcommand{\doiorurl}{%
	\iffieldundef{url}
	{\iffieldundef{eprint}
		{}
		{http://arxiv.org/abs/\strfield{eprint}}}
	{\strfield{url}}%
}

\newcommand{\myhref}[1]{%
	\ifboolexpr{%
		test {\ifhyperref}
		and
		not test {\iftoggle{bbx:eprint}}
		and
		not test {\iftoggle{bbx:url}}
	}
	{\href{\doiorurl}{#1}}
	{#1}%
}

\usepackage[
	bookmarks=true,
	linktocpage=true,
	bookmarksnumbered=true,
	breaklinks=true,
	pdfstartview=FitH,
	hyperfigures=false,
	plainpages=false,
	naturalnames=true,
	colorlinks=true,
	pagebackref=false,
	pdfpagelabels]{hyperref}

\hypersetup{
	colorlinks,
	citecolor=blue,
	filecolor=blue,
	linkcolor=blue,
	urlcolor=blue
}

\usepackage[capitalize, noabbrev]{cleveref}
\let\subsectionSymbol\S
\crefname{subsection}{\subsectionSymbol\!}{subsections}

\calclayout

\makeatletter
\@namedef{subjclassname@2020}{%
	\textup{2020} Mathematics Subject Classification}
\makeatother

\newtheorem{theorem}{Theorem}

\theoremstyle{definition}
\newtheorem{definition}[theorem]{Definition}
\newtheorem*{definition*}{Definition}

\newtheorem*{remark*}{Remark}
\newtheorem{example}[theorem]{Example}
\newtheorem*{example*}{Example}

\newtheorem*{convention*}{Convention}

\newtheorem*{notation*}{Notation}

\newtheorem*{question*}{Question}

\newcommand{\ot}{\otimes}

\newcommand{\Z}{\mathbb{Z}}

\renewcommand{\k}{\Bbbk}
\newcommand{\Sym}{\mathbb{S}}

\newcommand{\Ftwo}{{\mathbb{F}_2}}

\DeclarePairedDelimiter\set{\{}{\}}

\newcommand{\id}{\mathsf{id}}
\renewcommand{\th}{\mathrm{th}}

\newcommand{\Hom}{\mathrm{Hom}}

\newcommand{\xla}[1]{\xleftarrow{#1}}
\newcommand{\xra}[1]{\xrightarrow{#1}}
\newcommand{\defeq}{\stackrel{\mathrm{def}}{=}}

\newcommand{\bC}{\mathbb{C}}

\newcommand{\bF}{\mathbb{F}}

\newcommand{\cA}{\mathcal{A}}

\newcommand{\cH}{\mathcal{H}}

\newcommand{\cM}{\mathcal{M}}
\newcommand{\cN}{\mathcal{N}}

\newcommand{\rC}{\mathrm{C}}

\newcommand{\rK}{\mathrm{K}}
\newcommand{\rM}{\mathrm{M}}
\newcommand{\rP}{\mathrm{P}}
\newcommand{\rT}{\mathrm{T}}

\usepackage{caption, subcaption}
\usepackage{footmisc}
\usepackage[boxruled, noend]{algorithm2e}
\crefname{algocf}{Algorithm}{algocf}

\DeclareMathOperator*{\img}{\mathrm{img}}
\DeclareMathOperator{\rank}{rank}
\DeclareMathOperator{\card}{card}
\newcommand{\F}{\mathbb{F}_2}

\DeclareMathOperator{\xor}{\triangle}

\newcommand{\RP}{\mathbb{R}\mathrm{P}}

\newcommand{\Rplus}{\protect\hspace{-.1em}\protect\raisebox{.35ex}{\smaller{\smaller\textbf{+}}}}
\newcommand{\Cpp}{\mbox{C\Rplus\Rplus}\xspace}
\title{Persistence Steenrod modules}

\author[U.~Lupo]{Umberto Lupo}
\address{U.L., Institute of Bioengineering, School of Life Sciences, École Polytechnique Fédérale de Lausanne (EPFL), Lausanne, Switzerland}
\email{\href{mailto:umberto.lupo@epfl.ch}{umberto.lupo@epfl.ch}}

\author[A.~Medina-Mardones]{Anibal~M.~Medina-Mardones}
\address{A.M-M., Max Planck Institute for Mathematics \and University of Notre Dame}
\email{\href{mailto:ammedmar@mpim-bonn.mpg.de}{ammedmar@mpim-bonn.mpg.de}}

\author[G.~Tauzin]{Guillaume Tauzin}
\address{G.T., INAIT SA}
\email{\href{mailto:gtauzin@protonmail.com}{gtauzin@protonmail.com}}

\date{\today}
\subjclass[2020]{55N31, 55S10, 62R40, 68T09}
\keywords{Persistent homology, Steenrod squares, Steenrod barcode, cyclo-octane molecule}

\begin{document}

\begin{abstract}
    It has long been envisioned that the strength of the barcode invariant of filtered cellular complexes could be increased using cohomology operations.
    Leveraging recent advances in the computation of Steenrod squares, we introduce a new family of computable invariants on mod 2 persistent cohomology termed $Sq^k$-barcodes.
    We present a complete algorithmic pipeline for their computation and illustrate their real-world applicability using the space of conformations of the cyclo-octane molecule.
\end{abstract}
 	\maketitle

\section{Introduction}
\label{s:introduction}

Persistent homology is one of the primary tools in the rapidly developing field of topological data analysis.
A motivating example for this technique is the study of a finite point cloud of data embedded in Euclidean space.
From it, we can produce a collection of nested simplicial complexes
\begin{equation*}
	X_0 \to X_1 \to \cdots \to X_n.
\end{equation*}
For example, by taking the nerve of balls of uniformly increasing diameter whose centers are the given points.
The homology construction provides us with a collection of linear maps
\begin{equation} \label{e:persistent homology intro}
\begin{tikzcd}[column sep = small]
H_\bullet(X_0;\Bbbk) \arrow[r] & H_\bullet(X_{1};\Bbbk) \arrow[r] & \cdots \arrow[r] & H_\bullet(X_n;\Bbbk).
\end{tikzcd}
\end{equation}
This is an example of a (graded) persistence module and its barcode, a summary of the way Betti numbers are shared by consecutive simplicial complexes, serves as a principled and robust feature of the data.
This invariant is effectively computable and some of the open source software developed for this end can be found in \cite{medina2021giotto, gudhi, bauer2021ripser, pérez2021giotto-ph}.
For an expository treatment of persistent homology, we refer the reader to \cite{carlsson2005data} or \cite{computationaltopology2010}.

Assuming $\Bbbk$ is a field, a straightforward duality argument shows that the barcode of \eqref{e:persistent homology intro} is equivalent to the barcode of its persistent cohomology
\begin{equation} \label{e:persistent cohomology intro}
\begin{tikzcd}[column sep = small]
H^\bullet(X_0; \Bbbk) & \arrow[l] \cdots & \arrow[l] H^\bullet(X_{n-1}; \Bbbk) & \arrow[l] H^\bullet(X_n;\Bbbk).
\end{tikzcd}
\end{equation}

When $\Bbbk$ is the field $\mathbb F_p$ with $p$ elements, we can define for each $k \geq 1$ an additional barcode naturally associated to \eqref{e:persistent cohomology intro} by consistently recording the ranks of Steenrod's cohomology operation $P^k$ in the commutative diagram
\begin{equation*}
\begin{tikzcd}[column sep = small]
H^\bullet(X_0; \mathbb F_p) & \arrow[l] \cdots & \arrow[l] H^\bullet(X_{n-1}; \mathbb F_p) & \arrow[l] H^\bullet(X_n; \mathbb F_p) \\
H^\bullet(X_0; \mathbb F_p) \arrow[u, "P^k"] & \arrow[l] \cdots & \arrow[l] H^\bullet(X_{n-1}; \mathbb F_p) \arrow[u, "P^k"] & \arrow[l] H^\bullet(X_n; \mathbb F_p) \arrow[u, "P^k"].
\end{tikzcd}
\end{equation*}

In this work we focus on the case $p=2$ where Steenrod operations are denoted $Sq^k$ and referred to as Steenrod squares,  a term that comes from the fact that $Sq^k([\alpha]) = [\alpha] \smallsmile [\alpha]$ where $k$ is the cohomological degree of $\alpha$ and $\smallsmile$ denotes the cup product on cohomology.
The ranks of Steenrod operations and in particular of Steenrod squares are able to detect finer information beyond the Betti numbers of a space.
For example:
\begin{enumerate}
	\item The real projective plane and the wedge of a circle and a sphere have, with $\F$-coefficients, the same Betti numbers, yet the rank of $Sq^1$ tells them apart.
	\item Similarly, the complex projective plane and the wedge of a 2-sphere and a 4-sphere have the same Betti numbers with any coefficients, yet the rank of $Sq^2$ distinguishes them.
	\item The suspensions of the two spaces above have the same Betti numbers and also isomorphic cohomology rings, yet the rank of $Sq^2$ tells them apart.
\end{enumerate}

The main contribution of this work is the theoretical and algorithmic developments needed to use the finer discriminatory power of Steenrod squares in persistent cohomology.
Specifically, we define $Sq^k$-barcodes and introduce a method for their computation.
Using a performance oriented implementation of our methods, we present examples showing that the finer information they reveal is non-trivially present in the point cloud of conformations of the cyclo-octane molecule $\text{C}_{8}\text{H}_{16}$.

\subsection*{Outline}

We begin in \cref{s:preliminaries} with an overview of the basic notions used in the rest of this article.
They concern homological algebra, simplicial topology, and persistence theory.
In \cref{s:steenrod_squares} we introduce the Steenrod squares $Sq^k$ through explicit formulas at the cochain level and discuss their relevance.
In \cref{s:persistence_steenrod_modules} we introduce persistence Steenrod modules and their $Sq^k$-barcodes, and show how to computationally incorporate these invariants into the persistent cohomology pipeline.
We present three examples in \cref{ss:example}, including one obtained from data sampled from the space of conformations of $\text{C}_{8}\text{H}_{16}$.
We close by providing conclusions in \cref{s:conclusion}.

\subsection*{Acknowledgments}

We would like to thank Ulrich Bauer, Prasit Bhattacharya, Matteo Caorsi, Kathryn Hess, Maximilian Schmahl, Gard Spreemann, Dennis Sullivan, and Ulrike Tillmann for their suggestions, questions, and comments about this project.
We thank the reviewers for their careful reading and many keen suggestions.
We are grateful to Ingrid Membrillo Solis for providing useful datasets and insight on the topic of molecular conformational spaces.
We thank the hospitality of the \textit{Laboratory for Topology and Neuroscience} at EPFL and acknowledge partial support from Innosuisse grant \mbox{32875.1 IP-ICT-1}.

\section{Conventions and preliminaries} \label{s:preliminaries}

We assume familiarity with the notions of \textit{chain complex} over a ring $\k$ and of its associated \textit{homology} graded $\k$-module.

\subsection{Tensor and hom complexes}

In this subsection we review two natural chain complexes associated to any pair of chain complexes $C$ and $C^\prime$.

The \textit{tensor product} $C \ot C^\prime$ is the chain complex whose degree-$n$ part is
\[
\left(C \ot C^\prime\right)_n = \bigoplus_{i+j=n} C_i \ot C^\prime_j,
\]
where $C_i \ot C^\prime_j$ is the tensor product of $\k$-modules, and whose boundary map is defined by
\[
\partial (v \ot w) = \partial v \ot w + (-1)^{|v|} v \ot \partial w.
\]

The \textit{hom complex} $\Hom(C, C^\prime)$ is the chain complex whose degree-$n$ part is the subset of linear maps between them that increase degree by $n$, i.e.,
\[
\Hom(C, C^\prime)_n = \{f \colon C \to C^\prime \ |\ \forall k \in \Z, \ f(C_k) \subset C^\prime_{k+n}\},
\]
and boundary map defined by
\[
\partial f = \partial_{C^\prime} \circ f - (-1)^{n} f \circ \partial_C.
\]
A \textit{chain map} is a $0$-cycle in this chain complex, and two chain maps are \textit{chain homotopy equivalent} if they are homologous cycles.
We extend this terminology and say that two maps $f, g \in \Hom(C, C^\prime)$ are \textit{homotopic} if their difference is nullhomologous, referring to a map $h \in \Hom(C, C^\prime)$ such that $\partial h = f - g$ as a \textit{homotopy} between them.

Regarding $\k$ as a chain complex with $0$-part equal to $\k$ and all other parts $0$, the \textit{linear dual} of a chain complex $C$ is the chain complex $\Hom(C, \k)$.
For historical reasons we will use cohomological grading for the dual of a chain complex, placing the dual of a chain in degree $n$ also in degree $n$ instead of $-n$ as would be more appropriate.

For any three chain complexes, there is a natural adjunction isomorphism:
\begin{equation} \label{e:adjuntion isomorphism}
\Hom(C \ot C^\prime, C^{\prime\prime}) \cong
\Hom(C, \Hom(C^\prime, C^{\prime\prime})).
\end{equation}

\subsection{Invariants and coinvariants}

Symmetries on chain complexes play an important role in this work.
Let $G$ be a finite group.
We will later focus solely on the symmetric group $\Sym_2$.
We denote by $\k[G]$ the \emph{group ring} of $G$, i.e., the free $\k$-module generated by $G$ together with the ring product defined by linearly extending the group structure on $G$.
We refer to a chain complex of $\k[G]$-modules as a chain complex \emph{with a $G$-action}.

To any chain complex $C$ with a $G$-action we naturally associate the following two chain complexes.
The subcomplex of \textit{invariant chains} of $C$, denoted by $C^G$, contains all elements $c \in C$ satisfying $g \cdot c = c$ for every $g \in G$.
The quotient complex of \textit{coinvariant chains} of $C$, denoted by $C_G$, is the chain complex obtained by identifying elements $c, c^\prime \in C$ if there exists $g \in G$ such that $c^\prime = g \cdot c$.

Let $C$ and $C^\prime$ be chain complexes and assume $C$ has a $G$-action.
The chain complex $\Hom(C, C^\prime)$ has a $G$-action induced from $(g \cdot f)(c) = f(g^{-1} \cdot c)$ and there is an isomorphism:
\begin{equation} \label{e:invariant coinvariant ajunction}
\Hom(C, C^\prime)^G \cong \Hom(C_G, C^\prime).
\end{equation}

\subsection{Simplicial complexes}

Simplicial complexes are used to combinatorially encode the topology of spaces and occur naturally on real-world data.

An \textit{abstract and ordered simplicial complex}, or a \textit{simplicial complex} for short, is a pair of sets $(V, X)$ where $V$ is a poset and the elements of $X$ are non-empty finite subsets of $V$, such that:
\begin{enumerate}
	\item The restriction of the partial order of $V$ to any element in $X$ defines a total order on it.
	\item For every $v$ in $V$, the singleton $\{v\}$ is in $X$.
	\item If $x$ is in $X$ and $y$ is a subset of $x$, then $y$ is in $X$.
\end{enumerate}
We abuse notation and denote the pair $(V, X)$ simply by $X$.

The elements of $X$ are called \textit{simplices} and the \textit{dimension} of a simplex is defined by subtracting $1$ from its cardinality.
Simplices of dimension $d$ are called $d$-simplices.
We abuse terminology and refer to the elements of $V$ and to their associated $0$-simplices both as \textit{vertices}.

A simplicial complex $Y$ is a \textit{subcomplex} of a simplicial complex $X$ if every simplex of $Y$ is a simplex of $X$.
In this case we say that $(X, Y)$ is a \textit{simplicial complex pair} and write $(X, Y) \subseteq (X^\prime, Y^\prime)$ if $X \subseteq X^\prime$ and $Y^\prime \subseteq Y$.

A \textit{filtered simplicial complex} is a simplicial complex $X$ together with subcomplexes
\[
\emptyset = X_{-\infty} \subseteq \cdots \subseteq X_i \subseteq X_{i+1} \subseteq \cdots \subseteq X_{+\infty} = X.
\]

\subsection{Simplicial cohomology}

Let $(X, Y)$ be a simplicial complex pair.
Denoting the subsets of \mbox{$n$-dimensional} simplices by $X_n \subset X$ and $Y_n \subset Y$, the chain complex $C_\bullet(X, Y; \k)$ of \textit{relative chains} of the pair $(X, Y)$ is defined as follows: Its degree-$n$ part is
\[
C_n(X, Y; \k) =
\frac{\k \big\{ X_n \big\}}
{\k \big\{ Y_n \big\}}
\]
i.e., the $\k$-module freely generated by the $n$-dimensional simplices in $X$ modulo those in $Y$, and its differential, referred to as \textit{boundary map}, is defined on basis elements by
\[
\begin{tikzcd}[column sep=normal, row sep=tiny, row sep=0
,/tikz/column 1/.append style={anchor=base east}
,/tikz/column 2/.append style={anchor=base west}
]
C_n(X, Y; \k) \arrow[r, "\partial_n"] & C_{n-1}(X, Y; \k) \\
x\ \arrow[r, |->] & \sum_{i=0}^{n} (-1)^i d_i x
\end{tikzcd}
\]
where $d_i$ is the operator that removes the $i^\th$ element in $x$ with respect to the induced total order.
We refer to $C_\bullet(X, \emptyset; \k)$ simply as the \textit{absolute chains} of $X$ and use the notation $C_\bullet(X; \k)$.

The \textit{relative cochains} of the pair $(X, Y)$ is the cochain complex $C^\bullet(X, Y; \k)$ defined explicitly by
\[
C^{n}(X, Y; \k) = \Hom_{\k}(C_n(X, Y; \k),\ \k)
\]
with
\[
\delta_{n}(\alpha)(c) = (-1)^n\alpha(\partial_n c).
\]
We refer to $C^\bullet(X, \emptyset; \k)$ as the \textit{absolute cochains} of $X$ and use the notation $C^\bullet(X; \k)$ for it.
Notice that $C^{n}(X, Y; \k)$ is isomorphic to the subspace of $C^{n}(X; \k)$ that vanish on $C_n(Y; \k)$.
The \textit{cohomology} of this pair, denoted by $H^{\bullet}(X, Y; \k)$, is defines as the cohomology of $C^{\bullet}(X, Y; \k)$.

Given pairs $(X_1, Y_1) \subseteq (X_2, Y_2)$ there are natural maps:
\begin{align*}
C^{\bullet}(X_1, Y_1; \k) &\leftarrow C^{\bullet}(X_2, Y_2; \k),
\\
H^{\bullet}(X_1, Y_1; \k) &\leftarrow H^{\bullet}(X_2, Y_2; \k),
\end{align*}
respectively defined and induced by restriction.

Relative and absolute homology in the simplicial context is defined similarly but we do not use them in this work.

\subsection{Persistence theory}

In this subsection $\k$ is assumed to be a field.
We will now review the basic concepts of the theory of persistence over $\k$ from a point of view that prioritizes persistent cohomology.
We refer to \cite{carlsson2005barcode} or \cite{desilva2011duality} for a more detailed exposition.

The totally ordered set $\overline{\Z}$, known as \textit{extended integers}, is the union of $\Z$ with two elements $-\infty$ and $+\infty$ such that
\[
-\infty < i < +\infty
\]
for any integer $i$.

A \textit{persistence module} $\mathcal M$ (over $\k$) is a diagram of $\k$ vector spaces and linear maps
\[
\begin{tikzcd}[column sep = small]
\cM(-\infty) & \cdots \arrow[l] & \cM(i) \arrow[l] & \cM(i+1) \arrow[l] & \cdots \arrow[l] & \arrow[l] \cM(+\infty).
\end{tikzcd}
\]
For $i \leq j$ in $\overline{\Z}$ we denote by $\mathcal M_{i,j}$ the unique composition $\cM(i) \leftarrow \cM(j)$ in the diagram.
We say that $\cM$ is \textit{pointwise finite-dimensional} (\textit{p.f.d.}) if the dimension of $\cM(i)$ is finite for each $i \in \overline{\Z}$.

A \textit{graded persistence module} $\cM^\bullet = \{\cM^d\}_{d \in \Z}$ is a collection of persistence modules indexed by the integers.
We say $\cM^\bullet$ is p.f.d. if each $\cM^d$ is.

A \textit{morphism} of persistence modules is a diagram of vector spaces and linear maps
\[
\begin{tikzcd}[column sep = small, row sep = small]
\cN(-\infty) & \cdots \arrow[l] & \cN(i) \arrow[l] & \cN(i+1) \arrow[l] & \cdots \arrow[l] & \arrow[l] \cN(+\infty) \\
\cM(-\infty) \arrow[u] & \cdots \arrow[l] & \cM(i) \arrow[l] \arrow[u] & \cM(i+1) \arrow[l] \arrow[u] & \cdots \arrow[l] & \arrow[l] \cM(+\infty) \arrow[u].
\end{tikzcd}
\]
To any morphism $\phi$ of persistence modules we can naturally associate persistence modules corresponding to its kernel $\ker \phi$ and image $\img \phi$.

A \textit{multiset} is a pair $(M, \mu)$ where $M$ is a set and $\mu \colon M \to \overline \Z$ is a function attaining only values greater than $0$.
We refer to $\mu(m)$ as the \textit{multiplicity} of $m$ and define the \textit{cardinality} of a multiset $(M, \mu)$ as
\[
\card M = \sum_{m \in M} \mu(m)
\]
if this sum is defined and $+\infty$ otherwise.
We sometimes regard sets as multisets with multiplicity function constant and equal to $1$.

Let $\cM$ be a p.f.d. persistence module, its \textit{barcode} is the multiset $Bar_\cM$ of pairs $[p.q] \in \overline{\Z} \times \overline{\Z}$ such that for any two extended integers $i \leq j$
\[
\rank \cM_{i,j} = \card \big\{[p.q] \in Bar_{\cM} \mid p \leq i \leq j \leq q \big\}.
\]
The barcode is a complete invariant of p.f.d. persistence modules.
We will sometimes use the notation $(p-1, q] \defeq [p,q]$.
The \textit{finite} and \textit{infinite parts} of the barcode are defined by
\begin{align*}
Bar^{\, \mathrm{fin}}_{\cM} & = \set[\big]{[p.q] \in Bar_{\cM} \mid p,q \in \Z}, \\
Bar^{\, \mathrm{inf}}_{\cM} & = Bar_{\cM} \setminus Bar^{\, \mathrm{fin}}_{\cM}.
\end{align*}
The \textit{barcode} of a graded persistence module $\cM^\bullet$ is the collection
\[
Bar_{\cM^\bullet} = \big\{ Bar_{\cM^d} \big\}_{d \in \Z}\ .
\]

Given a filtered simplicial complex $X$, define respectively its \textit{persistent relative} and \textit{absolute cohomology} by
\begin{align*}
\mathcal H_R^\bullet(X; \k)(i) & = H^\bullet(X, X_{i}; \k), &
\mathcal H_A^\bullet(X; \k)(i) & = H^\bullet(X_{i}; \k),
\end{align*}
with linear maps induced by restriction.
When $X$ and $\k$ are clear from the context we omit them from the notation.
We say $X$ is p.f.d. if either, and therefore both, of these are.

The barcodes of these persistence modules contain equivalent information.
More precisely, as shown in \cite{desilva2011duality} or more categorically in \cite{bauer2020structure}, the finite parts of these are equal as graded multisets after a degree shift and there is a bijection of multisets between their infinite parts.
Explicitly,
\begin{equation} \label{e:duality of barcodes}
\begin{tikzcd}[column sep = small, row sep=-5]
Bar_{\mathcal H_R^n}^{\mathrm{fin}} \arrow[r, "\cong"] &
Bar_{\mathcal H_A^{n-1}}^{\mathrm{fin}} \\
{[p,q]} \arrow[r, |->] & {[p,q]}
\end{tikzcd}
\quad \text{and} \quad \
\begin{tikzcd}[column sep = small, row sep=-5]
Bar_{\mathcal H_R^n}^{\mathrm{inf}} \arrow[r, "\cong"] &
Bar_{\mathcal H_A^{n}}^{\mathrm{inf}} \\
{[-\infty, r]} \arrow[r, |->] & {[r, +\infty]}.
\end{tikzcd}
\end{equation}

Although we do not use them in this work, persistent relative and absolute homology can be defined similarly and shown to have barcodes containing the same information encoded by those above.

\section{Steenrod squares} \label{s:steenrod_squares}

In this section we introduce the cohomology operations
\[
Sq^k \colon H^\bullet(X, Y; \F) \to H^{\bullet}(X, Y; \F)
\]
defined for any simplicial complexes pair $(X, Y)$ and every integer $k$.
These operations are natural.
In particular, for $(X_1, Y_1) \subseteq (X_2, Y_2)$ the diagram
\begin{equation} \label{e:naturality of Steenrod squares}
\begin{tikzcd}
H^\bullet(X_2, Y_2; \F) \arrow[r, "Sq^k"] \arrow[d] & H^{\bullet}(X_2, Y_2; \F) \arrow[d] \\
H^\bullet(X_1, Y_1; \F) \arrow[r, "Sq^k"] & H^{\bullet}(X_1, Y_1; \F)
\end{tikzcd}
\end{equation}
commutes.
Therefore, as will be developed in \cref{s:persistence_steenrod_modules}, any $Sq^k$ defines an endomorphism of the persistent (absolute and relative) cohomology of a filtered complex.

\subsection{History and definition} \label{ss:history}

The diagonal map of spaces
\[
\begin{tikzcd}[row sep = 0]
X \arrow[r, "D"] & X \times X \\
x \arrow[r, |->]& (x, x)
\end{tikzcd}
\]
induces a product in cohomology with field coefficients
\[
\smallsmile \colon H^\bullet(X) \ot H^\bullet(X) \xra{\cong}
H^\bullet(X \times X) \xra{H^\bullet(D)}
H^\bullet(X),
\]
which is (graded) commutative, since the diagonal is invariant under the transposition
\[
\begin{tikzcd}[row sep = 0]
X \times X \arrow[r, "T"] & X \times X \\
(x, y) \arrow[r, |->]& (y, x).
\end{tikzcd}
\]
One can then ask if this product be defined with integer coefficients.
During the mid 1930's Alexander, Kolmogorov, \v{C}ech and Whitney \cite{whitney1988history} defined the cup product dualizing a simplicial chain approximation to $D$ given by
\[
\begin{tikzcd}[row sep = 0,
/tikz/column 1/.append style={anchor=base east},
/tikz/column 2/.append style={anchor=base west}
]
C_\bullet \arrow[r, "\Delta"] &
C_\bullet \ot C_\bullet \\
{[0, \dots, n]} \arrow[r, |->] &
\sum_{i=0}^{n} {[0, \dots, i] \ot [i, \dots, n]}.
\end{tikzcd}
\]
The chain map $\Delta$ is not invariant under the transposition map
\[
\begin{tikzcd}[row sep = tiny]
C_\bullet \ot C_\bullet \arrow[r, "T"] &
C_\bullet \ot C_\bullet \\
a \ot b \arrow[r, |->] &
(-1)^{|a||b|} b \ot a,
\end{tikzcd}
\]
that is to say, $\Delta - T \Delta \neq 0$.

In 1947, Steenrod published his seminal paper \cite{steenrod1947products} introducing the square operations through an effective construction of ``coherent homotopies" correcting the broken symmetry of $\Delta$ (denoted by $\Delta_0$ from now on).
To explain this, let us consider the map $(1 - T) \Delta_0$ as a $0$-cycle in $\Hom \left( C_\bullet, C_\bullet^{\ot 2} \right)$, a chain complex with an $\Sym_2$-action induced from $T$.
The \textit{cup-$1$ coproduct}, defined explicitly by
\[
\Delta_1 [0, \dots, n] =
\sum_{i<j} \pm \, [0, \dots, i, j, \dots, n] \ot [i, \dots, j],
\]
is a boundary for this cycle \big($\partial \Delta_1 = (1 - T) \Delta_0$\big).
The cup-$1$ coproduct $\Delta_1$ corrects the lack of symmetry of $\Delta_0$ homologically, but it is itself not symmetric.
Steenrod gave formulae for higher corrections, the \textit{cup-$i$ coproducts} $\Delta_i$, satisfying
\[
\partial (\Delta_{i+1}) = \Delta_i - (-1)^i T \Delta_i.
\]
More abstractly, if $W$ is the minimal resolution of $\Z$ by free $\Z[\Sym_2]$-modules
\[
\Z[\Sym_r]\{e_0\} \xla{1-T}
\Z[\Sym_r]\{e_1\} \xla{1+T}
\Z[\Sym_r]\{e_2\} \xla{1-T}
\cdots,
\]
he effectively constructed a natural equivariant chain map
\begin{equation} \label{e:steenrod diagonal}
W \ot C_\bullet \to C_\bullet^{\ot 2},
\end{equation}
where $C_\bullet$ denotes the chains of a simplicial complex.
Passing to mod $2$ coefficients, Steenrod extracted from this construction finer invariants on the cohomology of spaces which we now review.

Using the linear duality functor on the map \eqref{e:steenrod diagonal} and passing to invariant chains we have a chain map
\[
\Hom \left(C_\bullet \ot C_\bullet, \F \right)^{\Sym_2}
\longrightarrow
\Hom \left(W \ot C_\bullet, \F \right)^{\Sym_2},
\]
which we can complete, using the isomorphisms \eqref{e:adjuntion isomorphism} and \eqref{e:invariant coinvariant ajunction} of \cref{s:preliminaries}, to a commutative diagram
\[
\begin{tikzcd}
\Hom \left( C_\bullet \ot C_\bullet, \F \right)^{\Sym_2} \arrow[r] &
\Hom \left( W \ot C_\bullet, \F \right)^{\Sym_2} \arrow[d] \\
\left( C^\bullet \ot C^\bullet\right)^{\Sym_2} \arrow[u, "cross\ product"]&
\Hom \left( W_{\Sym_2} \ot C_\bullet, \F \right) \arrow[d] \\
C^\bullet \arrow[u, "doubleing"] \arrow[r, dashed]&
\Hom \left( W_{\Sym_2}, C^\bullet \right),
\end{tikzcd}
\]
where the choice of coefficients ensures the \textit{doubling map} $\alpha \mapsto \alpha \ot \alpha$ is linear.
Using the adjunction isomorphism \eqref{e:adjuntion isomorphism}, the dashed arrow defines a linear map
\[
\begin{tikzcd}[row sep=0, column sep = tiny]
C^\bullet \ot W_{\Sym_2} \arrow[r] &[-10pt] C^\bullet \\
\alpha \ot e_i \arrow[r, |->] & (\alpha \ot \alpha)\Delta_i(-)
\end{tikzcd}
\]
descending to mod $2$ homology, and the \textit{Steenrod square} operations are defined by reindexing this map.
Explicitly,
\[
\begin{tikzcd}[row sep=0, column sep=tiny]
Sq^k \colon H^n \arrow[r] & H^{n+k} \\
\phantom{Sq^k \colon}{[\alpha]} \arrow[r, |->] & \big[ (\alpha \ot \alpha)\Delta_{n-k}(-) \big].
\end{tikzcd}
\]
The importance of Steenrod operations in stable homotopy theory is hard to overstate, see for example \cite{adams1974stablehomotopy}.
For a more leisure exposition of the construction and properties of Steenrod squares we refer to, for example, \cite{tangora1968operations}.

\begin{remark*}
	The name of these operations comes from the fact that $Sq^k([\alpha]) = [\alpha] \smallsmile [\alpha]$ where $k$ is the cohomological degree of $\alpha$ and $\smallsmile$ denotes the cup product on cohomology.
	The non-triviality of Steenrod squares is an obstruction to the existence of a commutative product of cocycles lifting $\smallsmile$.
\end{remark*}

\begin{remark*}
	The operation $Sq^1$ agrees with the Bockstein homomorphism.
	Namely, the connecting homomorphism induced from the following exact sequence of coefficient
	\[
	0 \to \Ftwo \xra{2} \bF_4 \to \Ftwo \to 0.
	\]
\end{remark*}

\begin{remark*}
    Steenrod square operations are parameterized by classes on the mod~$2$ homology of $\Sym_2$.
    From this viewpoint, Steenrod defined operations at odd primes non-constructively using the mod $p$ homology of $\Sym_p$ \cite{steenrod1962cohomology}.
    We do not treat these operations in the present paper.
\end{remark*}

\subsection{Cup-$i$ formulas} \label{ss:definition}

Throughout the rest of this article we set the ground ring $\k$ to be the field with two elements $\F$.
We will describe explicitly a natural equivariant chain map
\[
W \ot C_\bullet \to C_\bullet \ot C_\bullet
\]
or, equivalently, an equivariant chain map
\[
\begin{tikzcd}[column sep=small, row sep=0,
/tikz/column 1/.append style={anchor=base east},
/tikz/column 2/.append style={anchor=base west}]
W \arrow[r] & \Hom\left(C_\bullet, C_\bullet \ot C_\bullet\right) \\
e_i \arrow[r, |->] & \Delta_i.
\end{tikzcd}
\]

Let $X$ be a simplicial complex and $x \in X_n$.
For a set
\[
U = \{u_1 < \dots < u_r\} \subseteq \{0, \dots, n\}
\]
we use the notation $d_U(x) = d_{u_1} \ldots\, d_{u_r}(x)$.

\begin{definition}[\cite{medina2021fast_sq}] \label{d:simplicial cup-i coproducts}
	The \textit{simplicial cup-$i$ coproduct}
	\[
	\Delta_i\colon C_\bullet(X; \F) \to C(X; \F)^{\ot 2}_\bullet
	\]
	is the linear map defined on a basis element $x$ in dimension $n$ by
	\begin{equation*}
	\Delta_i(x) = \sum_U d_{U^0}(x) \ot d_{U^1}(x),
	\end{equation*}
	where the sum is taken over all sets $U = \{u_1 < \dots < u_{n-i}\}$ with $u_j \in \{0, \dots, n\}$ and
	\[
	U^0 = \{u_j\ |\ u_j + j \equiv 0 \text{ mod } 2\}, \qquad
	U^1 = \{u_j\ |\ u_j + j \equiv 1 \text{ mod } 2\}.
	\]
\end{definition}

These formulas are in a sense dual to Steenrod's original in \cite{steenrod1947products} but, as shown in \cite{medina2022axiomatic}, they are equivalent.
We have the homological relation
\begin{equation*}
\partial(\Delta_{i+1}) = (1+T) \Delta_i
\end{equation*}
for any integer $i$, and naturality for pairs $(X_1, Y_1) \subseteq (X_2, Y_2)$ making the diagram
\begin{equation} \label{e:naturality of cup-i coproducts}
\begin{tikzcd}
C_\bullet(X_2, Y_2; \F) \arrow[r, "\Delta_i"] & C_{\bullet}(X_2, Y_2; \F)^{\ot 2} \\
C_\bullet(X_1, Y_1; \F) \arrow[u] \arrow[r, "\Delta_i"] & C_{\bullet}(X_1, Y_1; \F)^{\ot 2} \arrow[u]
\end{tikzcd}
\end{equation}
commute.

\begin{definition} \label{d:steenrod squares}
	Let $(X, Y)$ be a pair of complexes.
    The \textit{$k$-th Steenrod square}
	\[
	Sq^k\colon H^\bullet(X,Y;\F) \to H^{\bullet}(X,Y;\F)
	\]
	is the linear map sending a class $[\alpha]$ represented by a cocycle $\alpha \in C^{n}(X,Y;\F)$ to the class represented by the cocycle whose value on $c \in C_{n+k}(X,Y;\F)$ is
	\[
	(\alpha \ot \alpha) \Delta_{n-k} (c).
	\]
\end{definition}

We notice that thanks to \eqref{e:naturality of cup-i coproducts}, the Steenrod square operations are natural for pairs $(X_1, Y_1) \subseteq (X_2, Y_2)$, i.e., diagram \eqref{e:naturality of Steenrod squares} commutes.

\begin{example}
	Let us consider the model of the real projective plane $\RP^2$ presented in \cref{f:rp2 with chosen cocycle} together with the cocycle $\alpha$, dual to
	\[
	a = [2,4]+[2,3]+[3,5]+[1,5]+[1,4],
	\]
	representing the generator of $H^1(\RP^2; \F) \cong \F$.
	According to \cref{d:steenrod squares}, the cocycle $(\alpha \ot \alpha) \Delta_0(-)$ represents the class $Sq^1\big( [\alpha] \big) \in H^2(\RP^2)$.
	Using \cref{d:simplicial cup-i coproducts} and bilinearity, we are looking for basis elements $[i, j] \ot [i', j']$ appearing in $a \ot a$
	with $j = i'$ and such that $[i,j,j'] \in \RP^2$.
	The cocycle $(\alpha \ot \alpha) \Delta_0(-)$ is given by adding together $[i, j, j']$ for each such basis element.
	In our case, out of 25 basis elements appearing in $a \ot a$ only $[2, 3] \ot [3, 5]$ contributes a non-zero term and, therefore, $Sq^1\big([\alpha]\big)$ is represented by the cocycle dual to $[2,3,5]$, and
	\[
	\rank \big(Sq^1\colon H^1(\mathrm \RP^2;\F) \to H^2(\mathrm \RP^2;\F) \big) = 1.
	\]
\end{example}

\begin{figure}
	\begin{tikzpicture}[scale=2.3]
\coordinate (A) at (30: 1);
\coordinate (B) at (90: 1);
\coordinate (C) at (150: 1);
\coordinate (D) at (210: 1);
\coordinate (E) at (270: 1);
\coordinate (F) at (330: 1);

\draw[] (A) node[right] {6} -- (B) node[above] {4} -- (C) node[left] {2} --
(D) node[left] {6} -- (E) node[below] {4} -- (F) node[right] {2} -- (A);

\draw[] (A) -- (C) -- (E) -- (A);

\coordinate (X) at (90: .5);
\coordinate (Y) at (210: .5);
\coordinate (Z) at (330: .5);

\draw[] (X) -- (Y) -- (Z) -- (X);
\draw[] (X) -- (B) (Y) -- (D) (Z) -- (F);

\node at (90: .35) {3};
\node at (330: .35) {1};
\node at (210: .35) {5};

\draw[blue] (B) -- (C) -- (X) -- (Y) -- (Z) -- (E) -- (F);
\end{tikzpicture}
	\caption{Real projective plane $\RP^2$ together with a chosen representative of the non-zero class in $H^1(\RP2; \mathbb F_2)$.}
	\label{f:rp2 with chosen cocycle}
\end{figure}

\begin{remark*}
    The $\Delta_i$ maps are deeply rooted in the combinatorics of simplices.
    To illustrate their primitive nature we mentioned that another fundamental construction can derive from them: the nerve of $n$-categories \cite{street1987orientals, medina2020globular}.
    This is a reflection of a profound connection between convex geometry, higher category theory and Steenrod higher diagonals \cite{kapranov1991polycategory,medina2022fib_poly}.
\end{remark*}

\begin{remark*}
	We have focused on simplicial complexes since they are better known and lead to faster computations, but there are also effective constructions of Steenrod cup-$i$ coproducts for cubical complexes \cite{kadeishvili2003cupi, pilarczyk2016cubical, medina2021cubical}.
	Our algorithms, presented in \cref{ss:computations}, can be adapted using these to compute Steenrod barcodes of cubical complexes.
\end{remark*}

\begin{remark*}
    To define Steenrod operations effectively at any prime $p$, the \mbox{cup-$i$} coproducts where generalized in \cite{medina2021may_st} to cup-$(p,i)$ coproducts for simplicial and cubical chains using the operadic methods of P. May \cite{may1970general} and the model of the $E_\infty$-operad introduced by the second named author \cite{medina2020prop1, medina2021prop2}.
    These have been implemented in the computer algebra system \texttt{ComCH} \cite{medina2021comch}, and the incorporation of Steenrod operations at odd primes into the persistence pipeline is left to future work.
\end{remark*}

\subsection{Self-intersections} \label{ss:self-intersection}

\begin{figure}
	\newcommand*{\xMin}{0}%
\newcommand*{\xMax}{4}%
\newcommand*{\yMin}{0}%
\newcommand*{\yMax}{4}%
\begin{subfigure}{.4\textwidth}
	\centering
	\begin{tikzpicture}[scale=.8]
	\draw[-{Latex[length=2mm]}] (-.5,\yMin)--(-.5,\yMax);
	\draw[-{Latex[length=2mm]}] (-.5,\yMin)--(-.5,\yMax-.5);
	\draw[-{Latex[length=2mm]}] (4.5,\yMin)--(4.5,\yMax);
	\draw[-{Latex[length=2mm]}] (4.5,\yMin)--(4.5,\yMax-.5);
	
	\draw[-{Latex[length=2mm]}] (\xMin, -.5)--(\xMax, -.5);
	\draw[-{Latex[length=2mm]}] (\xMin, 4.5)--(\xMax, 4.5);
	
	\draw (0,0)--(0,4)--(4,4)--(4,0)--(0,0);

	\draw[color=blue!50, very thick] (0,2) .. controls (1,2.5) and (3,1.5) .. (4,2);
	\draw[color=red!50, very thick] (0,1) .. controls (1,1.3) and (3,1) .. (4,1);
	\end{tikzpicture}
	\caption{Torus}
\end{subfigure}
\quad  
\begin{subfigure}{.4\textwidth}
	\centering
	\begin{tikzpicture}[scale=.8]
	\draw[-{Latex[length=2mm]}] (-.5,\yMin)--(-.5,\yMax);
	\draw[-{Latex[length=2mm]}] (-.5,\yMin)--(-.5,\yMax-.5);
	\draw[-{Latex[length=2mm]}] (4.5,\yMax)--(4.5,\yMin);
	\draw[-{Latex[length=2mm]}] (4.5,\yMax)--(4.5,\yMin+.5);
	
	\draw[-{Latex[length=2mm]}] (\xMin, -.5)--(\xMax, -.5);
	\draw[-{Latex[length=2mm]}] (\xMin, 4.5)--(\xMax, 4.5);
	
	\draw (0,0)--(0,4)--(4,4)--(4,0)--(0,0);
		
	\draw[color=blue!50, very thick] (0,2) .. controls (1,2.5) and (3,1.5) .. (4,2);
	\draw[color=red!50, very thick] (0,1) .. controls (1,1) and (2,2.5) .. (4,3);
	\end{tikzpicture}
	\caption{Klein Bottle}
\end{subfigure}
	\caption{The torus is obtained by gluing together horizontal and vertical boundary components with no twists.
		The self-intersection for any 1-cycle is always even.
		The Klein bottle is obtained by gluing horizontal and vertical boundary components with a twist.
		The self-intersection of the depicted 1-cycle is odd.}
	\label{f:torus and klein bottle}
\end{figure}

From a geometric viewpoint, the cup product can be interpreted in terms of intersections of cycles in certain cases.
For any space, Thom showed that every mod $2$ homology class is represented by the push-forward of the fundamental class of a closed manifold $W$ along some map to the space.
Furthermore, if the target $M$ is a closed $n$-manifold, and therefore satisfies Poincar\'{e} duality
\[
H^k(M; \F) \xra{\cong} H_{n-k}(M; \F),
\]
The cohomology class dual to the homology class represented by the intersection of two transverse maps $V \to M$ and $W \to M$, or more precisely their pull-back $W \times_M V \to M$, is the cohomology class $[\alpha] \smallsmile [\beta]$ where $[\alpha]$ and $[\beta]$ are respectively dual to the homology classes represented by $W \to M$ and $V \to M$.

By taking $[\alpha] = [\beta]$ of cohomological degree $k$, we have that $Sq^k \big( [\alpha] \big) = [\alpha] \smile [\alpha]$ is represented by the transverse self-intersection of $W \to M$, that is, the intersection of this map and a generic perturbation of itself.

For example, let us consider the Torus $\rT$ and the Klein bottle $\rK$.
Two manifolds with the same mod~2 Betti numbers.
These surfaces are distinguished by the fact that
\begin{align*}
& \rank \big( Sq^1\colon H^1(\rT; \F) \to H^2(\rT; \F) \big) = 0, \\
& \rank \big( Sq^1\colon H^1(\rK; \F) \to H^2(\rK; \F) \big) = 1,
\end{align*}
which we can verify geometrically, as depicted in \cref{f:torus and klein bottle}, by noticing that the transverse self-intersection of a map $S^1 \to \rT$ parallel to any boundary component is always an even number of disjoint points, whereas that of a map $S^1 \to \rK$ parallel to the untwisted boundary component is always an odd number.

In manifold topology, the relationship at the (co)homology level between cup product and intersection is classical.
For a comparison between these at the level of (co)chain see \cite{medina2021flowing}.

\subsection{Relations and further structure}

The cup-$i$ coproducts or, equivalently, their linear dual products, arise from effectively constructing coboundaries that coherently enforce the commutativity relation of the cup product in cohomology.
This is an example of a general principle: constructing cochains enforcing cohomological relations lead to further cohomological structures.
In our case, the commutativity relation of cup product gives rise to Steenrod operations.

There are two notable relations satisfied by the Steenrod square operations.
The first one, known as the \textit{Cartan relation}, expresses the interaction between these operations and the cup product:
\[
Sq^k \big( [\alpha] [\beta] \big) =
\sum_{i+j=k} Sq^i \big( [\alpha] \big)\, Sq^j \big( [\beta] \big),
\]
whereas the second, the \textit{Adem relation} \cite{adem1952iteration}, expresses dependencies appearing among the iteration of operations:
\begin{equation} \label{e:adem relations}
	Sq^i Sq^j = \sum_{k=0}^{\lfloor i/2 \rfloor} \binom{j-k-1}{i-2k} Sq^{i+j-k} Sq^k
\end{equation}
where $\lfloor- \rfloor$ denotes the integer part function and the binomial coefficient is reduced mod $2$.

To tap into the secondary structure associated with these relations, one needs to provide effective proofs for them, that is to say, construct explicit cochains that enforce the relations when passing to cohomology.
Such effective proofs were recently given respectively in \cite{medina2020cartan} and \cite{medina2021adem}, and we expect that the additional structure they unlock will also play a role in applied topology.

\section{Persistence Steenrod modules} \label{s:persistence_steenrod_modules}

In this section we define persistence Steenrod modules and their associated barcodes.
We also introduce, for finite filtered simplicial complexes, a complete algorithmic pipeline for their computation.

\begin{definition} \label{d:persistence steenrod module}
	A \textit{persistence Steenrod module} (\textit{over} $\Ftwo$) is a graded persistence module $\cM^\bullet$ of $\F$ vector spaces together with a degree-$k$ endomorphism
	\[
	Sq^k \colon \cM^{\bullet} \to \cM^{\bullet}
	\]
	for each integer $k$, such that $Sq^k(i)$ satisfies the Adem relations \eqref{e:adem relations} for every $i \in \overline{\Z}$.
\end{definition}

Similar to how persistence modules can be thought of as modules over the polynomial algebra $\F[x]$, persistence Steenrod modules correspond to modules over the algebra $\cA[x]$ where $\cA$ is the \textit{Steenrod algebra}, the free algebra generated by symbols $Sq^k$ modulo the ideal generated by the Adem relations.

\subsection{Steenrod barcode}

The following computable invariant of persistence Steenrod modules is central to applications.

\begin{definition}
	For any integer $k$, the \textit{$Sq^k$-barcode} of a persistence Steenrod module $\cM^\bullet$, denoted by $Sq^kBar_{\cM^\bullet}$, is the barcode of the image persistence module $\img Sq^k$.
	We refer to the collection of all of these as the \textit{Steenrod barcode} of $\cM^\bullet$.
\end{definition}

The following example illustrates that, unlike barcodes of regular persistence modules, Steenrod barcodes of persistence Steenrod modules are not a complete invariant.
Let $M^\bullet$ be the graded vector space given by
\[
M^0 = \F\{x_0\}, \qquad M^1 = \F\{x_1\}, \qquad M^2 = \F\{x_2, y_2\},
\]
and equal to $0$ in all other degrees.
Let $\cM$ be the graded persistence module
\[
\cM^\bullet(i) =\begin{cases}
M^\bullet & i = 0, \\
0 & \text{otherwise.}
\end{cases}
\]
We make $\cM^\bullet$ into a persistence Steenrod module in two non-isomorphic ways, but with the same Steenrod barcodes, by defining
\[
Sq^2(x_0) = x_2, \qquad Sq^1(x_1) = x_2, \qquad Sq^0 = \id,
\]
and
\[
Sq^2(x_0) = x_2, \qquad Sq^1(x_1) = y_2, \qquad Sq^0 = \id.
\]
We thank Prasit Bhattacharya for suggesting this example.

The most prominent examples of persistence Steenrod modules are given by persistent relative and absolute cohomology of a filtered complex $X$.
In this case, denoting both $\mathcal H^\bullet_R$ and $\mathcal H^\bullet_A$ by $\mathcal H^\bullet$, we have
\[
Sq^0 Bar_{\mathcal H^\bullet(X;\, \F)} \cong
Bar_{\mathcal H^\bullet(X;\, \F)}
\]
since $Sq^0$ is the identity.

\subsection{Duality} \label{ss:duality}

The following example illustrates that, unlike the case of regular barcodes discussed in \eqref{e:duality of barcodes}, the Steenrod barcode of persistent relative and absolute cohomology need not completely determine each other.
Let $\rM$ be the M\"{o}bius band and consider the filtration $S^1 \to \rM$ where the circle is included as the boundary of $\rM$.
Given that $H^\bullet(\rM, S^1)$ is isomorphic to the reduced absolute cohomology of the real projective plane $\rM/S^1 \cong \RP^2$, and that $\rM$ is homotopy equivalent to its central circle, one can verify that the Steenrod barcode of absolute cohomology is empty but that of relative cohomology is not.

An important case where the Steenrod barcodes of persistent relative and absolute cohomology determine each other is when there are only finite bars in their regular barcodes.
More precisely, let $X$ be a p.f.d.\ filtered complex such that for some integer $n$ either
\[
Bar_{\mathcal H^n_R(X;\, \F)} =
Bar_{\mathcal H^n_R(X;\, \F)}^{\mathrm{fin}}
\]
or, equivalently,
\[
Bar_{\mathcal H^{n-1}_A(X;\, \F)} =
Bar_{\mathcal H^{n-1}_A(X;\, \F)}^{\mathrm{fin}} \,.
\]
Then, for every integer $k$ there is a bijection of multisets
\begin{equation} \label{e:bijection of finite steenrod bars}
\begin{tikzcd}[column sep = normal, row sep=-5
,/tikz/column 1/.append style={anchor=base east}
,/tikz/column 2/.append style={anchor=base west}
]
Sq^k Bar_{\mathcal H^n_R(X;\, \F)} \arrow[r, "\cong"] &
Sq^k Bar_{\mathcal H^{n-1}_A(X;\, \F)} \\
{[p,q]} \arrow[r, |->] & {[p,q]}.
\end{tikzcd}
\end{equation}

We illustrate the argument in an example that contains all the ideas of the proof.
For complete details we refer to the study of dualities in the categorical framework presented in \cite{bauer2020structure}.
Consider $Sq^1$ and a two stage filtration $X_1 \to X_2$ or, more explicitly,
\[
X_n = \begin{cases}
\emptyset & n < 1, \\
X_1 & n = 1, \\
X_2 & n > 1,
\end{cases}
\]
and the diagram
\[
\begin{tikzcd}
\xla{\delta_{n+1}} H^{n+1}(X_1) &
\arrow[l, "\, i_{n+1}^\ast"'] H^{n+1}(X_2) &
\arrow[l, "\, j_{n+1}^\ast"'] H^{n+1} (X_2, X_1) \xla{\delta_{n}} \\
\xla{\delta_{n}} H^{n}(X_1) \arrow[u, "Sq^1"'] &
\arrow[l, "i_n^\ast"'] H^{n}(X_2) \arrow[u, "Sq^1"] &
\arrow[l, "j_n^\ast"'] H^{n}(X_2, X_1) \arrow[u, "Sq^1"] \xla{\delta_{n-1}} \\
\xla{\delta_{n-1}} H^{n-1}(X_1) \arrow[u, "Sq^1"'] &
\arrow[l, "\ i_{n-1}^\ast"'] H^{n-1}(X_2) \arrow[u, "Sq^1"] &
\arrow[l, "\ j_{n-1}^\ast"'] H^{n-1}(X_2, X_1) \arrow[u, "Sq^1"] \xla{\delta_{n-2}}
\end{tikzcd}
\]
where the horizontal maps are part of the long exact sequence of the pair $(X_2, X_1)$.
Consider $\beta \in H^n(X_2,X_1)$ with $Sq^1 \beta \neq 0$.
Since all regular bars are finite $j^\ast_n \, \beta = 0$ so $j^\ast_{n+1} Sq^1 \beta = 0$ and we have an $Sq^1$-bar $(0, 1]_R^{n}$.
By exactness, $\delta_{n-1} \alpha = \beta$ for some $\alpha \in H^{n-1}(X_1)$ where $\delta_{n-1}$ is the $(n-1)^\th$ connecting homomorphism.
Since these commute with Steenrod squares, we have $Sq^1 \alpha \neq 0$.
Furthermore, $Sq^1 \alpha$ is not in the image of $i_n^\ast$ since otherwise $Sq^1 \beta$ would be $0$.
Therefore, there is a $Sq^1$-bar $(0,1]_A^{n-1}$.

Conversely, given $\alpha \in H^{n-1}(X_1)$ with $Sq^1 \alpha \neq 0$ the finiteness assumption implies that $\delta_n \, Sq^1 \alpha \neq 0$ so we have a $Sq^1$-bar $(0,1]_A^{n-1}$.
Denote by $\beta$ the element $\delta_{n-1} \alpha$ and notice that $Sq^1 \beta \neq 0$ with exactness implying $j_n^\ast Sq^1 \beta = 0$, so we have a $Sq^1$-bar $(0,1]_R^{n}$.

\subsection{Truncations} \label{ss:truncations}

Given a filtered complex $X$ and an integer $n$ there are two naturally associated filtered complexes $X_{\geq n}$ and $X_{\leq n}$ defined respectively by
\[
(X_{\geq n})_k =
\begin{cases}
X_n & k < n, \\
X_k & k \geq n,
\end{cases}
\qquad
(X_{\leq n})_k =
\begin{cases}
X_k & k \leq n, \\
X_n & k > n,
\end{cases}
\]
and referred to as the \textit{above} and \textit{below truncations} at $X_n$.
Persistent relative (resp.\ absolute) cohomology behaves well with respect to above (resp.\ below) truncations.
Explicitly, there exist canonical inclusions
\[
\begin{split}
Bar_{\cH^\bullet_R(X_{\geq n})} \to
Bar_{\cH^\bullet_R(X)}, \\
Bar_{\cH^\bullet_A(X_{\leq n})} \to
Bar_{\cH^\bullet_A(X)},
\end{split}
\]
and
\[
\begin{split}
Sq^\bullet Bar_{\cH^\bullet_R(X_{\geq n})} \to
Sq^\bullet Bar_{\cH^\bullet_R(X)}, \\
Sq^\bullet Bar_{\cH^\bullet_A(X_{\leq n})} \to
Sq^\bullet Bar_{\cH^\bullet_A(X)}.
\end{split}
\]
We remark that this form of ``stability'' of Steenrod barcodes may fail when considering persistent relative (resp.\ absolute) cohomology and  below (resp.\ above) truncations.
For example, consider the filtration $S^1 \to \rM \to \rC\rM$, where $\rC\rM$ is the cone on the M\"{o}bius band.
The Steenrod barcode of the relative absolute cohomology of this filtration is empty whereas, as discussed at the beginning of \cref{ss:duality}, its below truncation at $\rM$ is not.

\newpage
\subsection{Computing the Steenrod barcode} \label{ss:computations}

In this subsection we provide algorithms to compute the Steenrod barcode of the persistent relative cohomology of a finite filtered simplicial complex $X$
\[
\emptyset = X_{-1} \subset X_0 \subset X_1 \subset \cdots \subset X_m = X,
\]
together with a total order of its elements
\[
a_0 < a_1 < \cdots < a_m
\]
such that for all $j \in \set{0, \dots, m}$ we have
\[
X_j = \set{a_i \in X\ |\ i \leq j}.
\]
Most of this pipeline is applicable to other filtered cellular complexes, with the exception of \cref{a:stsq}.

\subsubsection{Regular barcode}

Let us begin by reviewing an effective construction of the barcode of the persistent relative cohomology of $X$.
Let $D$ be the matrix representing
\[
\partial \colon C_\bullet(X; \F) \to C_\bullet(X; \F)
\]
in the canonical ordered basis $\{a_0 < \cdots < a_m\}$.
We index columns and rows in this matrix starting at $0$, and denote $\overline j = m-j$ for all $j \in \set{0, \dots, m}$.
Consider $D^\perp$ defined by
\[
D^\perp_{p,\, q} = D_{\overline q,\, \overline p}.
\]
Notice that $D^\perp_{\leq j, \leq j}$ represents the coboundary of $C^\bullet(X, X_{\overline j-1}; \F)$.

\begin{figure}
	\input{aux/reduce}
	\caption{Column reduction algorithm}
	\label{f:reduce}
\end{figure}

Applying to $D^\perp$ a version of \cref{a:reduce} in \cref{f:reduce} that remembers the performed operations we produce a reduced matrix $R$ and an upper triangular invertible matrix $V$ satisfying
\[
R = D^\perp V.
\]
Denoting the $j$-th column of $R$ by $R_j$, let
\[
P = \set{j \mid R_j = 0}, \qquad
N = \set{j \mid R_j \neq 0} \qquad
E = P \setminus \set{\text{pivots of } R}.
\]
There exists a canonical bijection between the union of $N$ and $E$, and the barcode of persistent relative cohomology given by
\begin{align*}
	N \ni j &\mapsto \big[\, \overline j, \overline{\mathrm{pivot}\,R_j} \, \big] \in Bar_{\cH^{\dim(a_j)+1}_R}^{\mathrm{fin}} \\
	E \ni j &\mapsto \big[\! -1, \overline j \, \big] \in Bar_{\cH^{\dim(a_j)}_R}^{\mathrm{inf}}
\end{align*}
that provides a preferred cocycle representative for each of these bars:
\[
[i,j] \mapsto
\begin{cases}
	V_{\overline j}, & i = -1, \\
	R_{\overline i}, & i \neq -1.
\end{cases}
\]
More specifically, a basis for $H^\bullet(X, X_{\overline j-1})$ thought of as a subspace in the direct sum
\[
\ker \delta = \img \delta \oplus H^\bullet(X, X_{\overline j-1}; \Bbbk),
\]
is given by the set of cochains corresponding to the vectors in the union of
\[
\set[\big]{R_k \mid k \in N,\, j < \mathrm{pivot}(R_k)}
\quad \text{and} \quad
\set{V_i \mid i \in E,\, i \leq j},
\]
and a basis for $\img \delta$ is given by
\[
\set{R_i \mid i \in N,\, i \leq j}.
\]

\subsubsection{Steenrod barcode}

We now describe an effective construction of the Steenrod barcode of the persistent relative cohomology of $X$.
For any integer $k \geq 0$, let $\mathtt{sq^k}$ be an algorithm taking as input a vector corresponding to a cochain $\alpha \in C^n(X, X_i)$ and producing the vector corresponding to the cochain
\[
(\alpha \ot \alpha) \Delta_{n-k}(-).
\]
Such an algorithm, based on the explicit formulas of \cref{s:steenrod_squares}, is presented as \cref{a:stsq} in \cref{f:stsq}.
Let $Q^k$ be the square matrix with columns given by
\[
Q^k_i = \begin{cases}
\mathtt{sq^k}(V_i) & i \in E, \\
\mathtt{sq^k}(R_j) & i = \mathrm{pivot}(R_j), \\
0 & \text{otherwise}.
\end{cases}
\]

\begin{figure}
	\input{aux/stsq}
	\caption{Algorithm producing for a simplicial complex $X$, non-negative integer $n$, integer $k$ between $1$ and $n$, and cocycle $\alpha$, presented as a set $A \subseteq X_n$, a cocycle representing $Sq^k([\alpha])$ identified with a set $B \subseteq X_{n+k}$.
	We use the notation $S \xor S^\prime = S \cup S^\prime \setminus (S \cap S^\prime)$ and $index(S) = \set{index(s) \mid s \in S}$.}
	\label{f:stsq}
\end{figure}

For matrices $M$ and $N$ of dimensions $m \times p$ and $m \times q$ we define the $m \times (p+q)$ matrix $M \mid N$ by
\[
(M \mid N)_i = \begin{cases}
M_i, & i \leq p, \\
N_{i-p}, & i > p.
\end{cases}
\]

We now have all the elements needed to introduce \cref{a:st_bar} in \cref{f:st-barcode} whose output is the $Sq^k$-barcode of the persistent relative cohomology of $X$.
Intuitively, the step from $j-1$ to $j$ either adds a new non-zero coboundary $R_j$ (which implies $Q^k_j = 0$) or the image $Q^k_j$ of a persistent cocycle generator (which implies $R_j=0$).
In either case, we need to reduce with respect to the subspace of coboundaries, generated by $R_{\leq j}$, the image of $Sq^k$, which is generated by $Q^k_{\leq j}$.
This process is done keeping track of when columns in $Q^k$ become zero and extracting from this information the $Sq^k$-barcode of the filtration.

\begin{figure}
	\input{aux/st_bar}
	\caption{Algorithm producing the $Sq^k$-barcode of a filtered simplicial complex given its reduced anti-transposed boundary matrix $R$ and a matrix $Q^k$ containing as a columns the images under $Sq^k$ of cocycles representing the barcode of its persistent relative cohomology.}
	\label{f:st-barcode}
\end{figure}

We leave the development of a pipeline for persistent absolute cohomology to future work, remarking that, as described in \eqref{e:bijection of finite steenrod bars}, its associated Steenrod barcode is equal to that of persistent relative cohomology if all regular bars are finite, a situation often countered in practice.


\section{Examples} \label{ss:example}

To demonstrate the feasibility of extracting Steenrod barcodes from realistic datasets using the computational pipeline described in \cref{s:persistence_steenrod_modules}, we have produced two open-source software implementations:
one\footnote{\label{f:steenroder python}Available at \url{https://github.com/Steenroder/steenroder}.} is a \texttt{Python} package optimized by means of the \texttt{Numba} library \cite{siukwan2015numba}, and the other\footnote{Available at \url{https://github.com/Steenroder/steenroder_cpp}} is a performance-oriented \Cpp package inspired by the \texttt{PHAT} library \cite{bauer2014phat}.

While detailed performance benchmarking is beyond the scope of this paper, some remarks are in order.
First, we note that both our implementations apply the \emph{clearing} optimization \cite{chen2011twist} to \cref{a:reduce}.
It is well-known (see e.g.\ the discussion in \cite{bauer2021ripser}) that clearing is particularly effective when computing relative persistent cohomology -- and even more so when the filtration is constructed via a Vietoris--Rips process.
Second, although the computation of matrix $Q^{k}$ (input to \cref{a:st_bar}) is in principle \emph{embarrassingly parallelizable} by tasking fully independent threads with the calculation of different columns, we have not yet pursued this path in our code.

Third, we expect (and observe experimentally) that ``sparsifying'' our filtrations via \emph{simplicial collapses} \cite{pritam2020collapses, boissonnat2020edge}, operations that preserve the homotopy type of each complex, can lead to a cascade of space and time improvements across our computational pipeline.
This is presumably because:
\begin{enumerate}
    \item the run-time and memory usage in \cref{a:reduce} is reduced, yielding sparser $R$ and $V$ matrices and hence cocycle representatives with smaller sizes on average;
    \item \label{i:smaller cocycles} the outer \textbf{forall} loop in \cref{a:stsq} becomes faster for smaller cocycles;
    \item the leaner $R$ matrix reduces the computational run-time and memory usage once again in the final Steenrod barcode computation, \cref{a:st_bar}.
\end{enumerate}
\cref{i:smaller cocycles} above deserves more emphasis: our experiments suggest that, in typical datasets, the main bottleneck in the entire pipeline is the computation of $\mathtt{sq}^k$ for a few exceptionally sizable cocycle representatives.
This is due to the quadratic complexity of the \textbf{forall} loop in \cref{a:stsq}.
By replacing the largest cocycle representatives returned by (any implementation of) \cref{a:reduce} with cohomologous ones with a smaller size, one could presumably alleviate this problem.
Our preliminary attempts using right-to-left reductions on the matrix $R$ output by \cref{a:reduce} have yielded promising results; in the future, we hope to further improve our implementations in this direction, as well as making it easily accessible through its incorporation into \texttt{giotto-tda} \cite{medina2021giotto}.

We now report the results of computing $Sq^1$-barcodes in a synthetic and a natural dataset.
In both cases we start from a point cloud and construct an associated filtered simplicial complexes through the Vietoris--Rips process with a fixed simplex dimension threshold of 3 (simplices with 4 or less vertices) and some distance threshold.
We close this section with a comparison of the $Sq^2$-barcode of two filtered complex models of the cone on the suspension of, respectively, $\bC\rP^2$ and $S^2 \vee S^4$.
Our experiments are fully reproducible as \texttt{Jupyter} notebooks.\footref{f:steenroder python}

\subsection{Flat Klein bottle} \label{ss:flat Klein bottle}

Our first example is constructed from a matrix of geodesic distances among $N$ points in a \emph{metrically flat} Klein bottle.
This is the Riemannian manifold $\mathcal{M} = (\mathbb{R}^2 / {\sim}, g)$ obtained from $\mathbb{R}^2$ with its usual metric via the equivalence relation $(x, y) \sim (x + n, 1 - y + m) \ \forall \ m, n \in \Z$.
To define this point cloud we selected $N = 100$ points corresponding to the vertices of a square grid inside the unit square $[0, 1]^2$.

\subsubsection{Persistent relative cohomology}

\begin{figure}[h!]
	\centering
	\includegraphics[width=\textwidth]{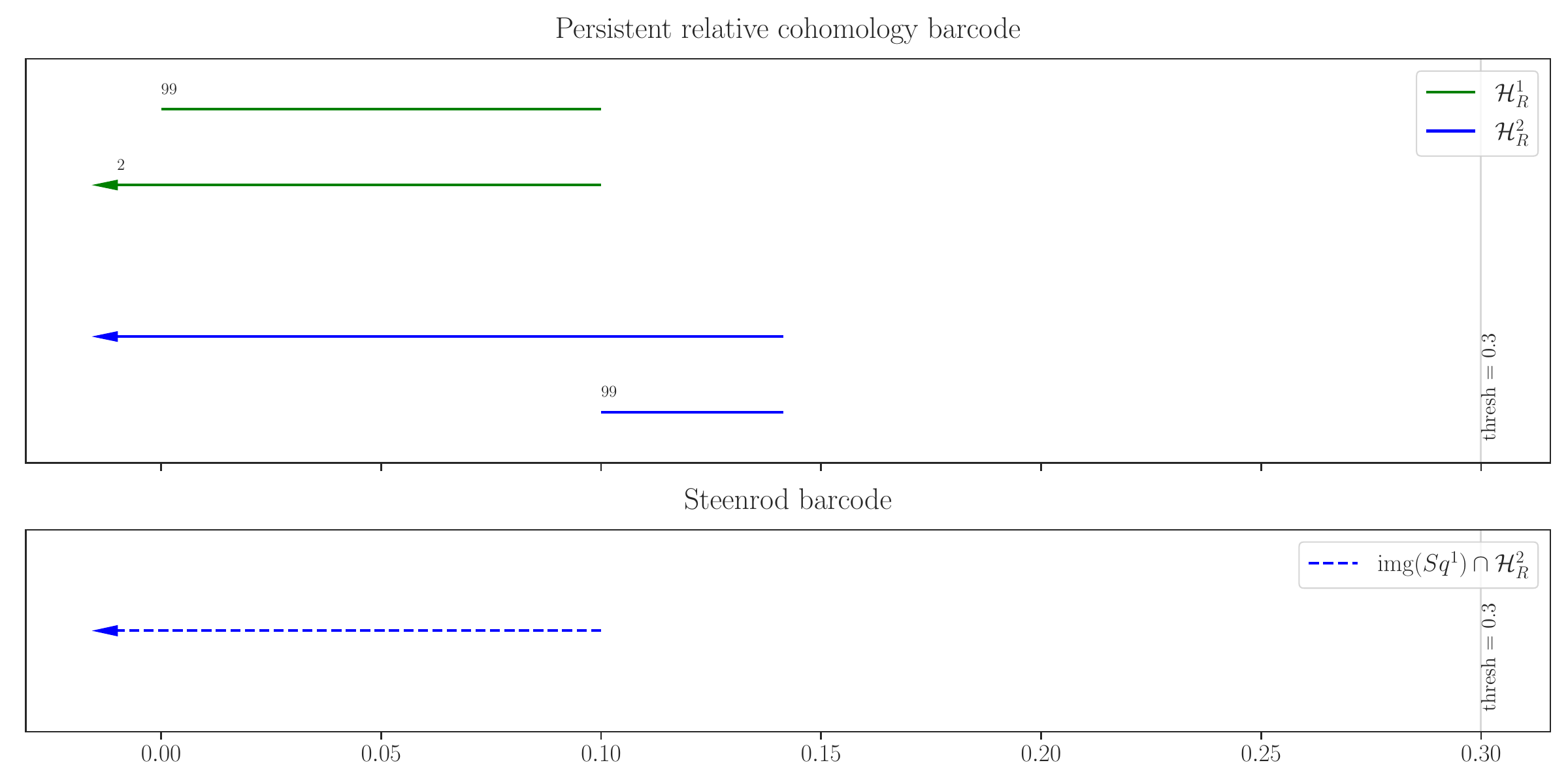}
	\caption{Regular barcode of $\mathcal H^i_R(X;\, \F)$ for $i = 1,2$ and associated $Sq^1$-barcode.
	Here $X$ is the Vietoris--Rips filtered complex associated to $N = 100$ points in a metrically flat Klein bottle and distance threshold $R = 0.3$.
	Arrowheads indicate infinite bars and integers over a bar denote its multiplicity (otherwise, the multiplicity is $1$).
    Vietoris--Rips filtration values are shown on the horizontal axis.}
	\label{f:flat Klein bottle relative}
\end{figure}

Let $X$ be the Vietoris--Rips filtered complex associated to this point cloud with distance threshold $R = 0.3$.
We apply our pipeline to compute the regular barcode of $\mathcal H^i_R(X;\, \F)$ for $i = 1,2$ and their associated $Sq^1$-barcode.
The results are presented in \cref{f:flat Klein bottle relative}.
There are three infinite bars: two in degree $1$ -- which happen to have identical birth and death due to the symmetry in our construction -- and one in degree $2$.
Our implementation detects an infinite bar in $\img(Sq^1) \cap \cH^2_R$.
This cohomological profile agrees with that expected from a filtered Klein bottle.

\subsubsection{Persistent absolute cohomology}

\begin{figure}
	\centering
	\includegraphics[width=\textwidth]{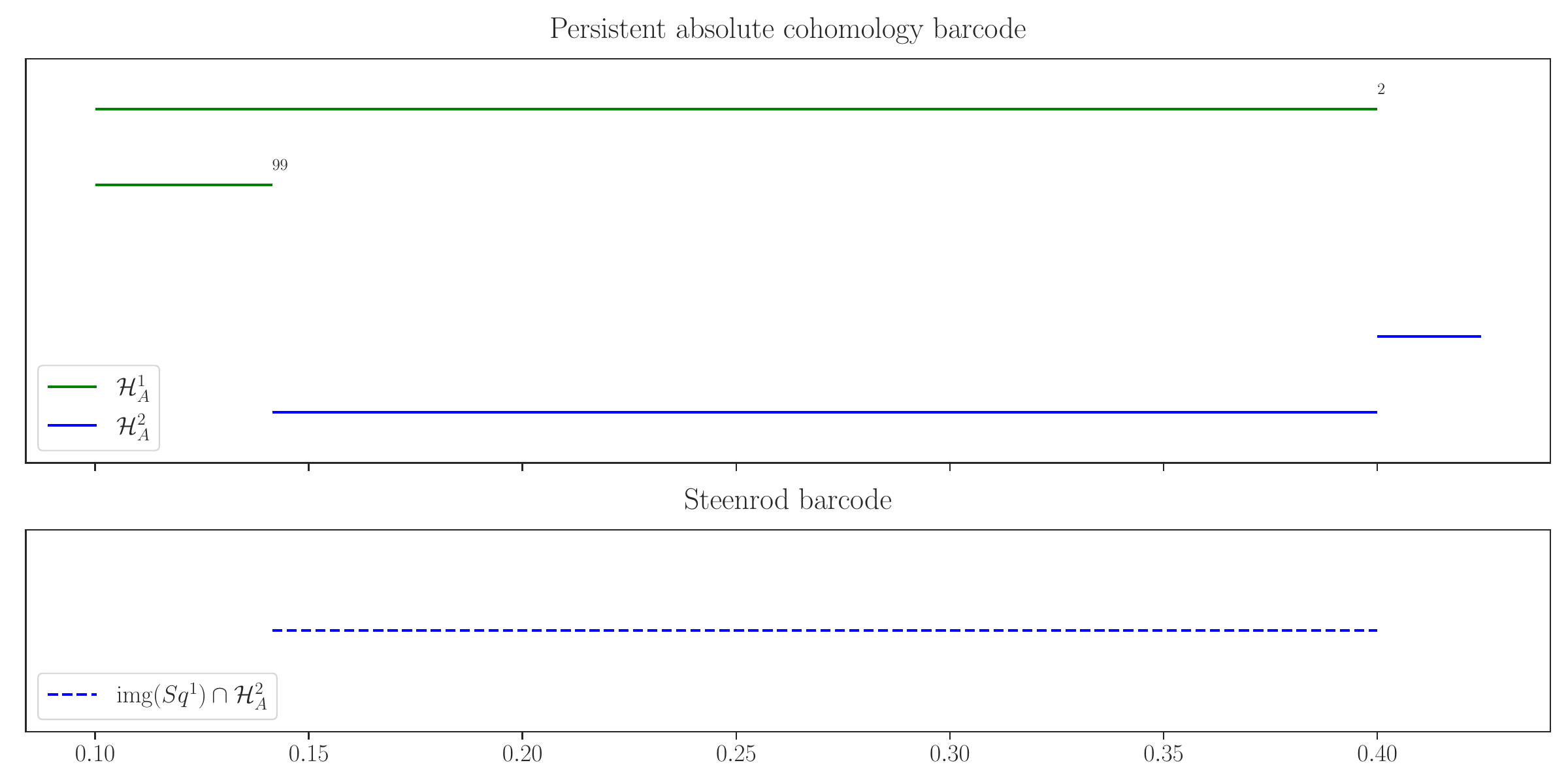}
	\caption{Regular barcode of $\mathcal H^i_A(X;\, \F)$ for $i = 1,2$ and associated $Sq^1$-barcode.
	Here $X$ is the Vietoris--Rips filtered complex associated to $N = 100$ points in a metrically flat Klein bottle and no distance threshold.}
	\label{f:flat Klein bottle absolute}
\end{figure}

Let $X$ be the Vietoris--Rips complex obtained from the same point cloud with no distance threshold.
In \cref{f:flat Klein bottle absolute} we present the regular barcode of $\mathcal H^i_A(X;\, \F)$ for $i = 1,2$ and associated $Sq^1$-barcode, obtained using our pipeline and the duality of Steenrod persistent cohomology (\cref{ss:duality}).
Our implementation detects a single Steenrod bar in $\img(Sq^1) \cap \cH^2_A$, which is born with the $\cH^1_A$ bars and dies with the $\cH^2_A$ bar.
Once again, this cohomological profile is consistent with that of a filtered Klein bottle.

\subsection{Conformational space of cyclo-octane}

Our second example involves a sampling of the conformational space of the cyclo-octane molecule $\text{C}_{8} \text{H}_{16}$.
We started with a dataset, originally from \cite{martin2010topology}, which consists of $6040$ vectors in $\mathbb{R}^{24}$.
Each of these vectors collects the 3D coordinates of all $8$ carbon atoms in a given cyclo-octane conformation after alignment to a reference one.
In \cite{martin2010topology}, this dataset was used to argue that the full conformational space of cyclo-octane is not a manifold, being in fact the union of a $2$-sphere with a Klein bottle glued together along two circles of singularities.
The reader can consult the papers \cite{membrillosolis2019topology, adams2021topology} for further details and references.

Candidate singular points in this dataset can be identified in a variety of ways; we used a set of $627$ singular points isolated in \cite{stolz2020geometry} via local persistent cohomology.\footnote{Data retrieved from \url{https://github.com/stolzbernadette/Geometric-Anomalies}.}
We removed these points from the dataset, and clustered the remaining $5413$ points using the HDBSCAN algorithm \cite{campello2013density} to obtain (samplings of) four $2$-strata -- presumably corresponding to a dense open subset of the Klein bottle, and three open connected subsets of the $2$-sphere.

As computed in \cite{membrillosolis2019topology}, the persistent absolute homology $\mathcal H_\bullet$ of the Vietoris--Rips filtered complex associated to the point cloud with $N = 3547$ elements supported on the presumed Klein bottle has, excluding the basic bar from $\mathcal H_0$, three prominent bars in its barcode.
Their birth and death values occur respectively before and after the value $R = 1.2$ and two come from $\mathcal H_1$ while the other from $\mathcal H_2$.

\subsubsection{Persistent relative cohomology}

\begin{figure}
	\centering
	\includegraphics[width=\textwidth]{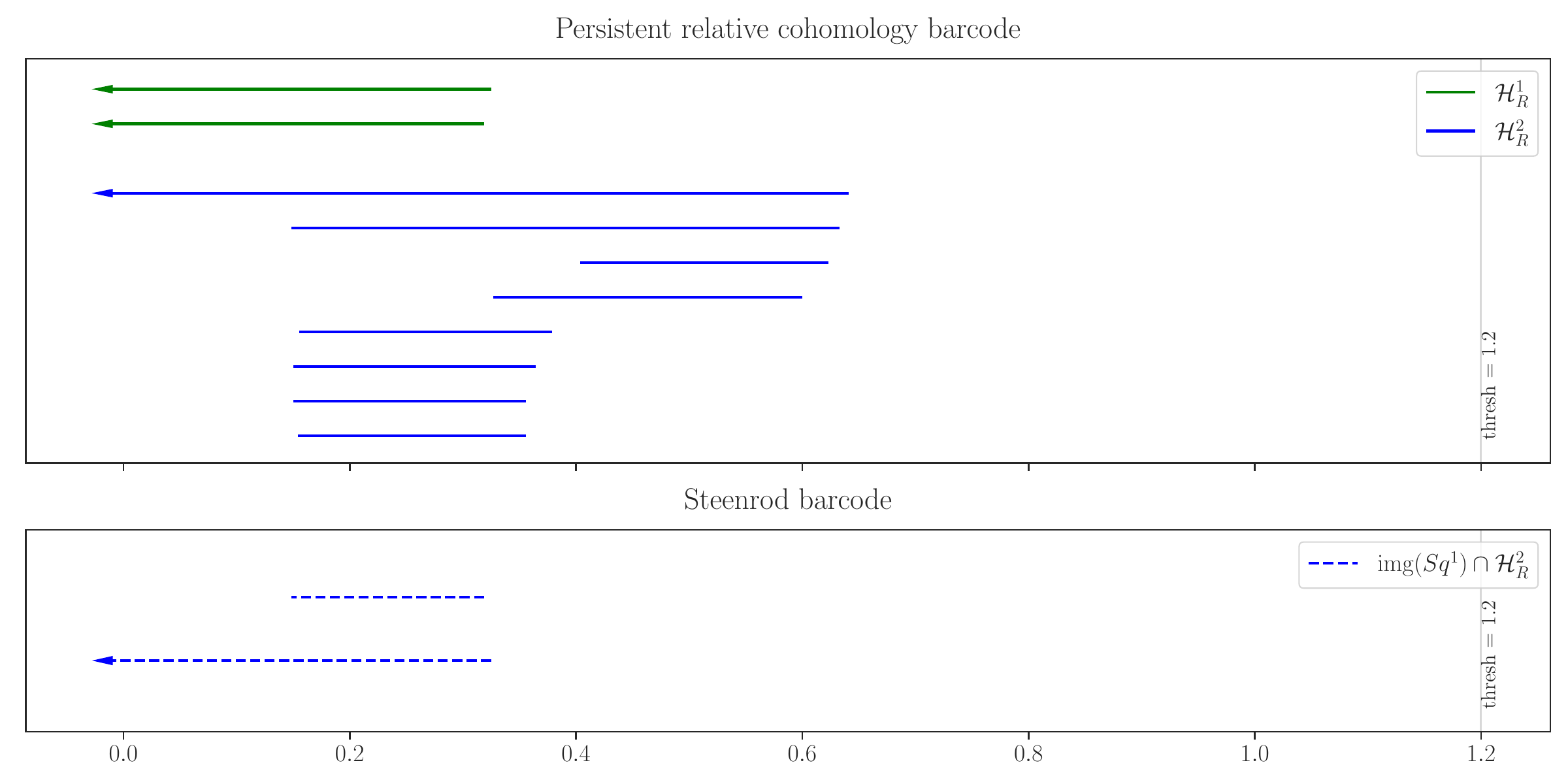}
	\caption{Regular barcode of $\mathcal H^i_R(X;\, \F)$ for $i = 1,2$ and associated $Sq^1$-barcode
		Here $X$ is the $R = 1.2$ Vietoris--Rips filtered complex associated to the $N = 3547$ ``Klein bottle component'' of the sampled conformational space of $\text{C}_{8} \text{H}_{16}$.
		Bars in the persistent relative cohomology barcode with lifetime shorter than $0.2$ are not shown to reduce clutter.}
	\label{fig:cyclo-octane}
\end{figure}

We study the persistent relative cohomology of $X$, the Vietoris--Rips filtered complex associated to this point cloud with distance thresholds $R = 1.2$.
In \cref{fig:cyclo-octane} we show, discarding short-lived bars ($< 0.2$) for ease of visualization, the regular barcode of $\mathcal H^i_R(X;\, \F)$ for $i = 1,2$ and their associated $Sq^1$-barcode.
Our implementation detects two $Sq^1$-bars.
One is infinite and born with one of the two infinite $\mathcal H^1_R$ bars, while the other is born with the other $\mathcal H^1_R$ bar and dies with the most prominent finite $\mathcal H^2_R$ bar.
The infinite parts of these barcodes are consistent with a filtered Klein bottle, where one of the infinite degree 1 bars interacts non-trivially with the degree 2 one.
The finite Steenrod bar adds extra information revealing a non-trivial interaction between the other infinite degree $1$ bar and a finite degree $2$ bar.
Refinements to the model for the conformation space of the cyclo-octane molecule resulting from the incorporation of this finer feature go beyond the scope of this example and are left unexplored.

\subsubsection{Persistent absolute cohomology}

\begin{figure}
	\centering
	\includegraphics[width=\textwidth]{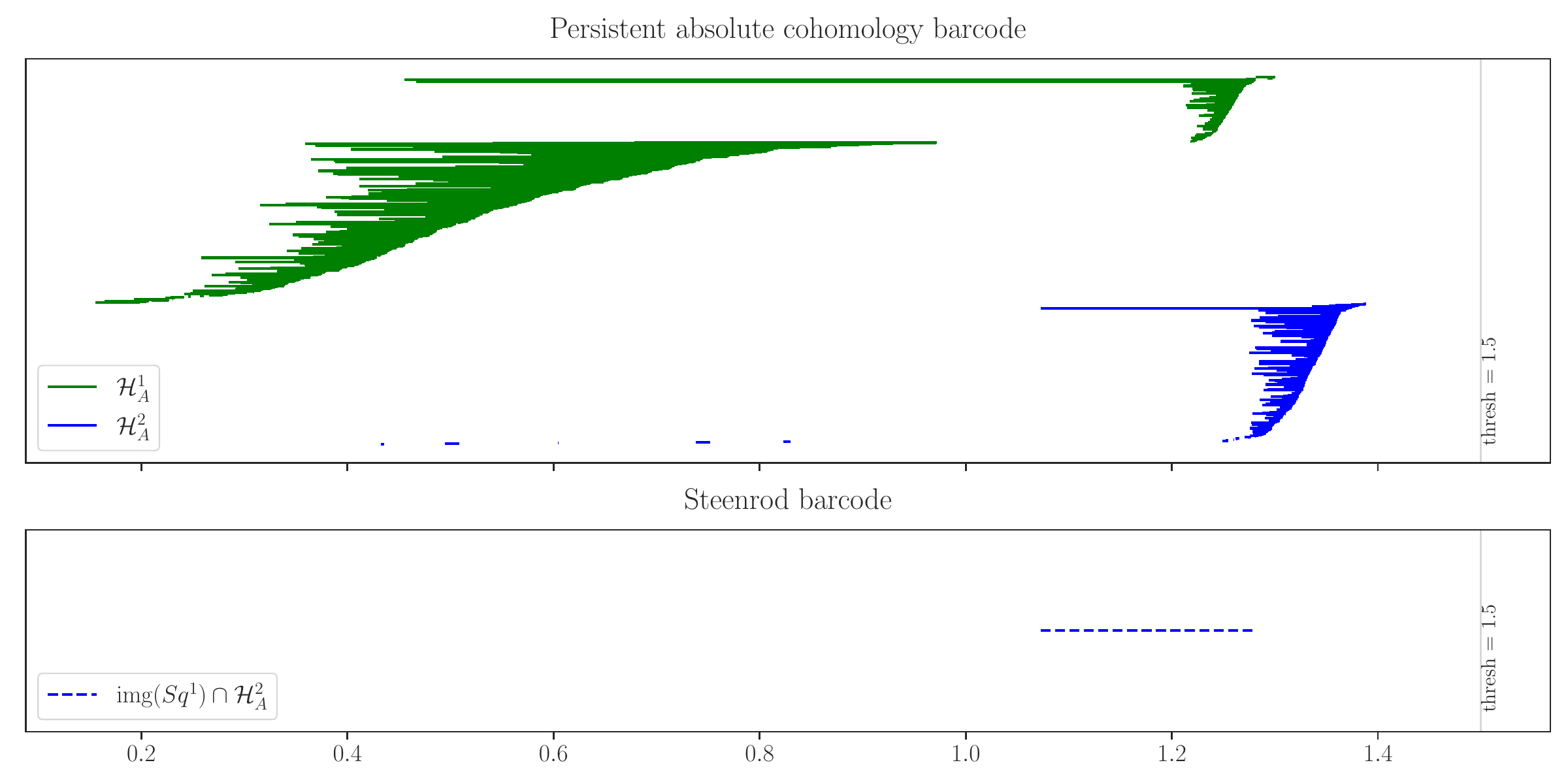}
	\caption{Regular barcode of $\mathcal H^i_A(X;\, \F)$ for $i = 1,2$ and associated $Sq^1$-barcode.
		Here $X$ denotes the $R = 1.5$ Vietoris--Rips filtered complex associated to an $N = 800$ random subsample of the ``Klein bottle component'' of the conformational space of $\text{C}_{8} \text{H}_{16}$.}
	\label{fig:cyclo-octane subsampled absolute}
\end{figure}

For this data set we will also compute a Steenrod barcode for absolute persistent cohomology.
We will use a distance threshold chosen to be larger than the death value of all prominent features.
To do so we consider a subsample consisting of $N = 800$ randomly selected points.
The persistent homology barcode of the associated Vietoris--Rips filtered complex with no distance threshold contains three prominent bars, and their death values are all less than $R = 1.5$.
As expected, two of these are associated to $\mathcal H_1$ and the other to $\mathcal H_2$.
Let $X$ be the the Vietoris--Rips filtered complex with distance thresholds $R = 1.5$ obtained from this subsample.
We remark that the threshold chosen ensures a correspondence between persistent relative and absolute cohomology of the Steenrod bars associated to their prominent features (\cref{ss:duality}).
In \cref{fig:cyclo-octane subsampled absolute} we show, discarding no bars, the regular barcode of $\mathcal H^i_A(X;\, \F)$ for $i = 1,2$ and associated $Sq^1$-barcode.
As expected, the interaction between the more prominent bars witnessed by the Steenrod barcode is consistent with a filtered Klein Bottle.

\subsection{Complex projective space and a wedge of spheres}

We conclude this section comparing the $Sq^2$-barcodes of the persistent absolute cohomology of two filtered simplicial complexes.
On one hand, we have a filtration of the cone on the suspension of the complex projective space $\mathrm C \, \Sigma \, \mathbb C\mathrm P^2$, and, on the other, one of $\mathrm C \, \Sigma \, (S^2 \vee S^4)$, where, as usual, $S^n$ denotes the $n$-dimensional sphere.
As mentioned in the introduction, $\Sigma \, (S^2 \vee S^4)$ and $\Sigma \, \mathbb C\mathrm P^2$ have isomorphic cohomology rings over any coefficients, but they can be distinguished by the action of $Sq^2$ on their mod~2 cohomology.

\begin{figure}
	\centering
	\begin{subfigure}[b]{0.49\textwidth}
		\centering
		\includegraphics[width=\textwidth]{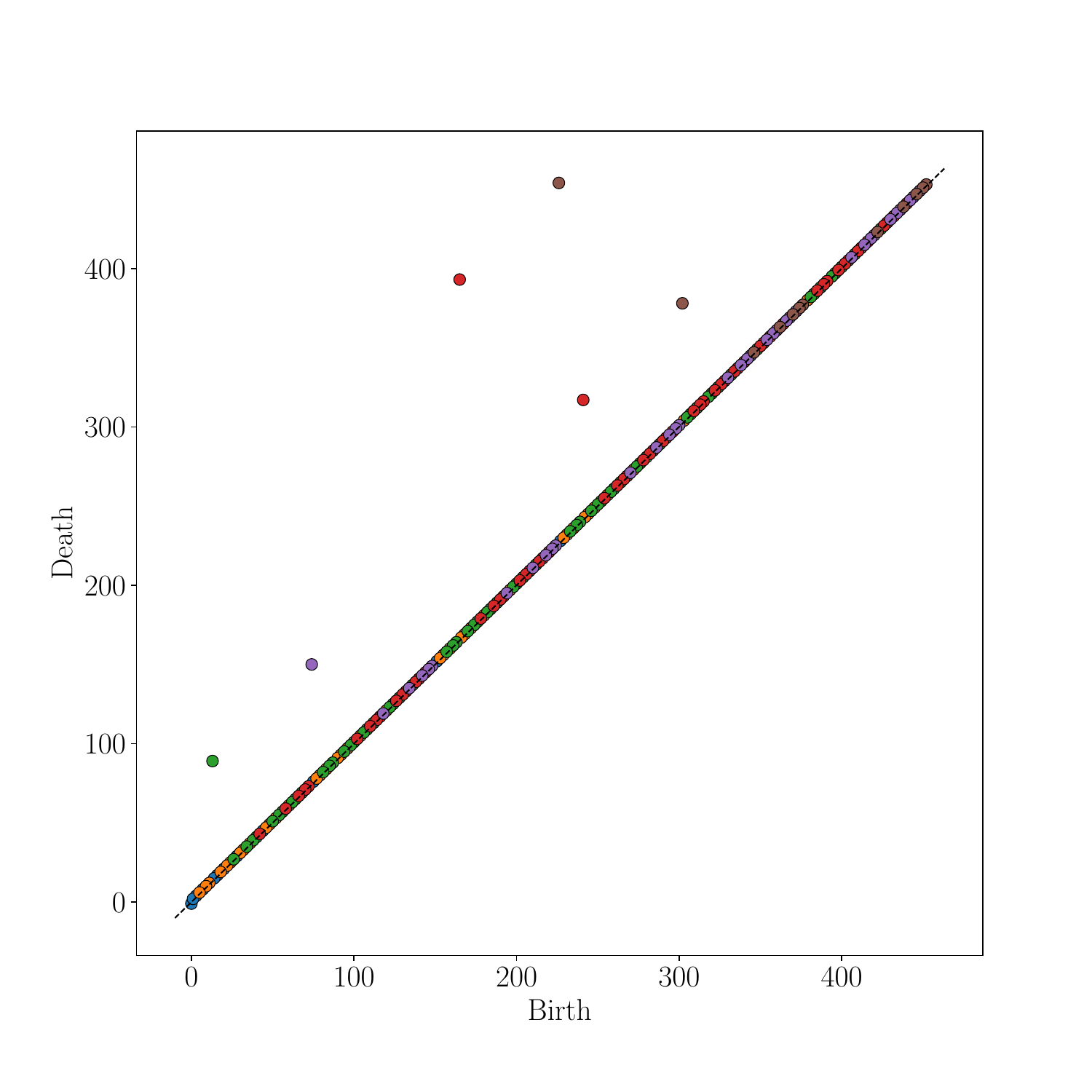}
		\caption{$\rC\,\Sigma(S^2 \vee S^4)$}
		\label{f:s2_s4}
	\end{subfigure}
	\begin{subfigure}[b]{0.49\textwidth}
		\centering
		\includegraphics[width=\textwidth]{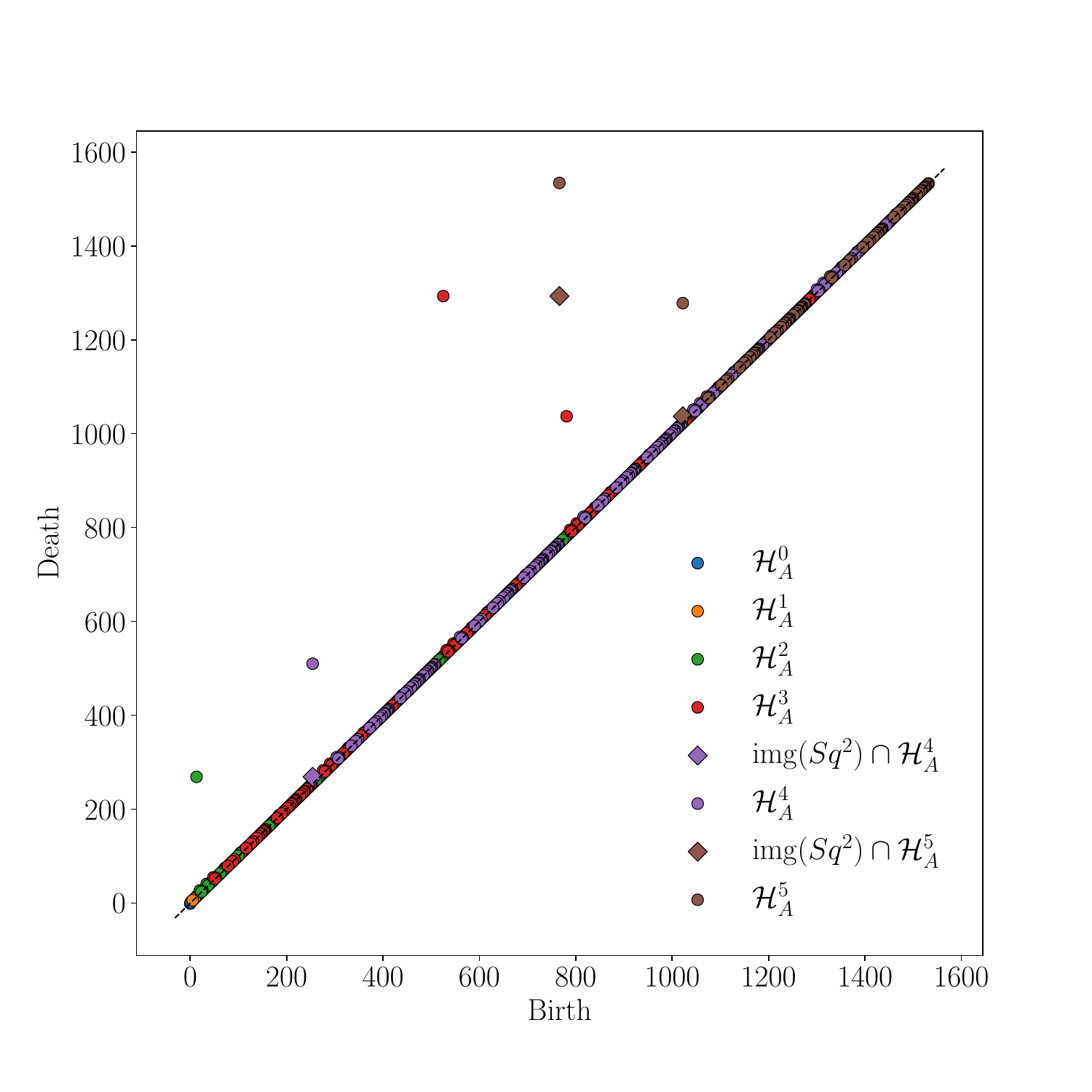}
		\caption{$\rC\,\Sigma\,\bC\rP^2$}
		\label{f:cp2}
	\end{subfigure}
	\caption{Persistence diagram representations of the regular and $Sq^2$-barcodes of the persistent absolute cohomology of (\textsc{a}) a filtered complex modeling the cone on the suspension of $S^2 \vee S^4$, and (\textsc{b}) a filtered complex modeling the cone on the suspension of $\bC\rP^2$.}
	\label{f:cp2 s_2s4}
\end{figure}

Interpreting bars as points in the plane, we plot in \cref{f:cp2 s_2s4} the regular and $Sq^2$-barcodes of the persistent absolute cohomology of these filtrations.
After rescaling by the number of simplices, we can see that the regular barcodes, symbolized by colored circles, are very similar;
yet there is a $Sq^2$-bar, represented by a brown diamond, present in the second figure only.


\section{Conclusion} \label{s:conclusion}

Steenrod barcodes increase the discriminatory power of traditional barcodes,
providing finer computable topological invariants of filtered spaces.
Furthermore, as we showed using the conformation space of $\text{C}_8\text{H}_{16}$, the additional information these invariants reveal is non-trivially present in real-world examples.
	\sloppy
	\printbibliography
\end{document}